\documentclass[final,12pt]{article}

\usepackage{tabularx}
\usepackage{amsfonts}
\usepackage{amsthm}
\usepackage{latexsym}
\usepackage{amsmath}
\usepackage{amssymb}
\usepackage{mathrsfs}
\usepackage{bm}

\def \bR {\mathbb R}
\def \bC {\mathbb C}

\def \al  {\alpha}

\def \n    {\mathfrak n}

\def \z    {\mathfrak z}
\def \v    {\mathfrak v}

\def \a    {\mathfrak a}
\def \s    {\mathfrak s}

\def \l    {\mathfrak l}

\def \lan  {\left\langle}
\def \ran  {\right\rangle}

\def \lb    {\left(}
\def \rb    {\right)}
\def \lq    {\left\lbrack}
\def \rq    {\right\rbrack}
\def \lg    {\left\lbrace}
\def \rg    {\right\rbrace}

\def \and {\quad \text{and} \quad}

\def \ha#1 {\hat{#1}}
\def \br#1#2 {\lq #1,#2\rq}
\def \Jl#1#2#3 {\lq \lq #1,#2\rq, #3 \rq}
\def \Jr#1#2#3 {\lq #1, \lq #2,#3\rq \rq}
\def \sc#1#2 {\lan #1,#2\ran}

\def \set#1 {\lg{#1}\rg}
\def \CO#1 {\com\lg{#1}\rg}

\def \bN {\mathbb N}

\newtheorem{theorem}{Theorem}[section]

\newtheorem{proposition}{Proposition}[section]

\title{Tangential Convergence of Bounded Harmonic Functions on Generalized Siegel Domains}

\author{M. Sundari\thanks{Chennai Mathematical Institute, Plot H-1,
SIPCOT IT Park, Padur P O, Siruseri, Chennai 603 103. Primary subject classification:22E, Secondary subject classification: 43A20}    }

\vspace{1.5cm}
\date{{\it To my parents}}

\usepackage{amsmath}
\begin{document}
\maketitle

\begin{abstract}
Suppose $u(x,y)$ is a bounded harmonic function on the upper half
plane such that $\lim\limits_{x\rightarrow\infty}u(x,y_o)=a$ for
some $y_o > 0$. Then one can prove that $\lim\limits
_{x\rightarrow\infty}u(x,y)=a$ for any other positive $y$. In this
paper, we shall consider the algebra of radial integrable
functions on the H-type groups and obtain a similar result for the
bounded harmonic functions on generalized Siegel domain.
\end{abstract}


\section{Introduction} \vskip0.15in \hskip0.15in Suppose
$u(x,y)$ is a bounded function on the upper half plane satisfying
$\Delta u=0$ and $\lim\limits_{x\rightarrow\infty}u(x,y_o)=a$ for
some positive number $y_o$. Then using classical methods, we can
prove that $\lim\limits _{x\rightarrow\infty}u(x,y)=a$ for any
other positive $y$. Here $\Delta$ denotes the Laplacian in the two
variables $x$ and $y$.

Suppose $H_r$ denotes the Heisenberg group of homogeneous dimension
$2r+2$. It is known that $H_r$ acts on the Siegel domain $D_r =
\{(z,z_o) \in \bC^r \times \bC : Im~~z_o
> \mid z\mid^2 \}$ by
translations. Under this action, $H_r$ is identified with the
boundary of $D_r$. By abuse of notation, we shall denote the
Laplace-Beltrami operator for the Bergman metric on $D_r$ by $\Delta$. Then we have
the following: \\
\vskip0.15in \begin{theorem}\label{1.1} Let $u$ be a bounded
function on $D_r$ such that $\Delta u=0$. If for an $\epsilon_o
> 0$, $\lim\limits_{(z,t)\rightarrow \infty} u(z,t,\epsilon_o)=a$,
then for any $\epsilon >0$ we have $\lim\limits_{(z,t)\rightarrow
\infty} u(z,t,\epsilon)=a$.
\end{theorem} \vskip0.15in
For any unexplained notation and terminology, one can refer to
\cite{HR}. The proof in the case of Heisenberg group depends on the
explicit form of the Poisson kernel and the Gelfand spectrum of the
commutative Banach algebra $L^1(H_r)^\natural$, of integrable radial
functions on $H_r$. A similar result for rank one symmetric spaces has 
been proved by Cygan \cite{JC}. Here we shall prove an analogue of Theorem
\ref{1.1} for H-type groups.

The paper is organized as follows: in section 2, we shall define a
H-type group and collect all the facts we require about H-type
groups. In section 3,we shall define the generalized Siegel domains
and describe the action of H-type groups on them. In section 4, we
shall prove our main result as a consequence of a Tauberian theorem
for H-type groups. The method of proof is that of
Hulanicki and Ricci \cite{HR}.

In the coming sections, we shall use the `variable constant
convention' according which our constants are denoted by $C$, $C'$
etc and these are not necessarily equal at different occurrences.

The author would like to thank Michael Cowling for his encouragement and many useful
discussions on the subject of this paper. This work has been done during her visit
to The University of New South Wales in June 2005. \vskip0.25in
\section{H-type groups}\vskip0.15in
\hskip0.15in In this section, we shall collect all the necessary information about
the H-type groups, N, and describe the Gelfand
transform for biradial functions on $N$. We shall denote the sets
$\bR \setminus \{0\}$, $\{0,1,2,\cdots\}$ by $\bR^*$, $\bN$
respectively and the semi-infinite interval $(0,\infty)$ by $\bR^+$.
For more detailed information on the material covered in this
section, one may refer to \cite{AB1} and the references there of.

Let $\n$ be a real two-step nilpotent Lie algebra endowed with an
inner product $\lan \cdot,\cdot \ran_{\n}$, let $\z$ be the centre
of $\n$. Write $\n$ as an orthogonal direct sum of two subspaces
$\v$ and $\z$ that is, $\n = \v \oplus \z$. For each $Z \in \z$,
define the map $J_Z:\v \rightarrow \v$ by the formula
\begin{align}
\lan J_Z X,Y\ran_\n = \lan [X,Y],Z\ran_\n,~~~ X,Y \in \v.
\end{align}
The Lie algebra $\n$ is said to be H-type if for every $Z \in \z$,
\begin{align}\label{e1.2}J_Z^2 = - \mid Z\mid^2
I_{\v},\end{align}
where $I_{\v}$ denotes the identity
transformation on $\v$. A connected, simply connected Lie group $N$
whose Lie algebra is H-type is said to be a H-type group. By
(\ref{e1.2}), we can see that every unit element $Z$ in $\z$ induces a complex
structure on $\v$ via the map $J_Z$. Therefore, $\v$ has even
dimension, say $2m$. If $k$ denotes the dimension of the centre $\z$
of $N$, then $Q=m+k$ is the homogeneous dimension of $N$.

As $N$ is a connected, simply connected nilpotent group, we know that
the exponential map $\exp : \n \rightarrow N$ is surjective.
Therefore, we shall identify $N$ with $\v\oplus\z$ and denote a
typical element $n$ of $N$ by $(X,Z)$ where $X \in \v$, $Z \in
\z$. Using the Campbell-Baker-Hausdorff formula, we get the product
rule in $N$ as
$$(X,Z)(X_1,Z_1)=(X+X_1,Z+Z_1+\displaystyle{\frac{1}{2}[X,X_1]}),
~~~~ X,X_1 \in \v, Z,Z_1 \in \z.$$ If $dX$ and $dZ$ denote the
Lebesgue measures on $\v$ and $\z$ respectively, then $dn=dX dZ$
denotes a Haar measure on $N$.

There are two classes of irreducible unitary representations of a
H-type group. Some are trivial on the centre and factor into
characters on $\v$. The others are parametrized by $\bR^+ \times
S_\z$ (-see \cite{CH} and \cite{AB1}) where $S_\z$ denotes the unit
sphere in $\z$. For $w$ in $S_\z$, we consider $\v$ endowed with
complex structure $J_w$. Denote by $I_w : \v \rightarrow \bC^m$ the
corresponding isomorphism. Then the corresponding Hermitian inner
product is given by
$$\{X, X_1\}_w = <X,X_1>_{\n}+ i < J_w X, X_1>_{\n},~~~X,X_1 \in \v.$$
Define
\begin{align*}
{\mathscr H}_{\nu,w}=\lbrace \xi : \v \rightarrow \bC\, |\,\xi \circ
I_w^{-1} : \bC^m \rightarrow \bC ~ \mbox{is
entire},\,
\|\xi\|_{\nu}^2 = \int\limits_{\v} \mid \xi(X)\mid^2 e^{-\frac{\nu
\mid X\mid^2}{2}}\, dX < \infty \rbrace
\end{align*}
Thus ${\mathscr H}_{\nu,w}$ is a Hilbert space with respect to the
inner product associated with the norm $\|\cdot\|_{\nu}$. For any
multi-index $j$ in $\bN^m$, we denote by ${\mathcal P}_{\nu, j}$ the
following normalized polynomial:
$${\mathcal P}_{\nu, j}(X)= \pi^{-m/2} (
\displaystyle{\frac{\nu}{2}})^{(m+\mid
j\mid)/2}(j!)^{-1/2}(I_w(X))^j,~~~X \in \v,$$ where $\mid j\mid
=j_1+ \ldots + j_m$, $j! = j_1! \ldots j_m!$ and $\zeta^j =
\zeta_1^{j_1} \ldots \zeta_m^{j_m}$ for $\zeta$ in $\bC^m$. One can
check that the family $\{{\mathcal P}_{\nu, j}\}_{j \in \bN^m}$ is an
orthonormal basis of ${\mathscr H}_{\nu,w}$.

For any $\nu$ in $\bR^+$ and any $w$ in $S_\z$, let $\pi_{\nu, w}$
be the unitary representation of $N$ on ${\mathscr H}_{\nu,w}$
defined by
$$[\pi_{\nu,w}(X,Z) \xi](X_1) = \displaystyle{e^{-\nu ( \frac{\mid X\mid^2}{4} + \frac
{\{X_1,X\}_w}{2} + i <Z,w>_{\n})}} \xi(X+X_1),~~X_1 \in \v,~\xi \in {\mathscr
H}_{\nu,w}.$$ Given $f
\in L^1(N)$, we shall define the group Fourier transform of $f$ as an
operator valued function on ${\mathscr H}_{\nu,w}$, by
$$\pi_{\nu,w}(f) = \int\limits_N \pi_{\nu,w}(n) f(n)\, dn.$$

A function $f$ on $N$ is said to be biradial if $f$ is radial in
both the variables $X$ and $Z$. In other words, there exist a
function $f_0$ on $\bR^2$ such that $f(X,Z)=f_0(\mid X\mid, \mid
Z\mid)$, $(X,Z) \in N$. Let $L^1(N)^\natural$ be the space of all
biradial integrable functions on the group $N$. We know from
\cite{DR} that $L^1(N)^\natural$ is a commutative Banach algebra.
The Gelfand spectrum of this commutative algebra is well known (see
\cite{K}, \cite{DR} and \cite{AB}) and can be described as follows:\\
\vskip0.1in
Let ${\mathscr J}_z$ be the generalized Bessel function defined for
every $x$ in $\bR$ by the rule
\begin{align*}
{\mathscr J}_z(x)= \begin{cases}
\displaystyle{\frac{\Gamma(z+1)}{\Gamma(\frac{2z+1}{2})\Gamma
(\frac{1}{2})}}\int\limits_{-1}^1 e^{ixs} (1-s^2)^{(2z-1)/2} ds
&\text{if}~~z > -1/2\ , \\ \cos x &\text{if}~~z=-1/2\ .
\end{cases}
\end{align*}
and let $L^\al_l$ be the $l^{th}$ Laguerre polynomial of order $\al$,
that is, $$L^\al_l(x)=\sum\limits_{j=0}^l \lb\begin{array}{c} l+\al\\
l-j\end{array}\rb \displaystyle{\frac{(-x)^j}{j!}},~~~~ x \in
\bR.$$ The bounded spherical functions of the commutative algebra
$L^1(N)^\natural$ are given by
\begin{align*} \phi_{\nu,l}(X,Z)&=&\displaystyle{e^{-
\frac{1}{4}\nu\mid X\mid^2}}
\displaystyle{\frac{L^{m-1}_l(\frac{1}{2}\nu\mid X\mid^2)}{\lb
\begin{smallmatrix}l+m-1\\l\end{smallmatrix}\rb}{\mathscr J}_{\frac{k-2} {2}}(\nu\mid
Z\mid)},~~~~~(X,Z) \in N,\\\phi_\mu(X,Z)&=&{\mathscr
J}_{m-1}(\mu \mid X\mid),~~~~~~~~~~~~~~~~~~~~~~~~(X,Z) \in N, \end{align*}
where $\nu >0$, $\mu \geq 0$ and $l \in \bN$. If $f$ is a biradial
integrable function on $N$, we have the Gelfand transform $\hat{f}$
of $f$ as a function on $\bR^+ \times \bN$ defined by the
rule:\begin{align}\label{1.4} \hat{f}(\nu,l)=\int\limits_N
f(n)\phi_{\nu,l}(n) \,dn~~~~\nu > 0,~~ l \in \bN.\end{align}
\vskip0.25in
\section{Harmonic
NA spaces and Poisson kernel} \vskip0.15in \hspace{1.5cm} In this
section, we shall describe the harmonic  $NA$ spaces, define
harmonic functions on the generalized Siegel domain and give the
explicit form of the Poisson kernel. For any unexplained terminology
and notation in this section, the reader may refer to \cite{DR} and
\cite{CDKR}. Let $S = NA$ be the semidirect product of the groups
$N$ and $A=\bR^+$ with respect to the action of $A$ on $N$ given by
the dilations $\delta_a :(X,Z) \mapsto (a^{1/2}X, aZ)$. We shall denote the Lie
algebras of 
$A$, $S$ by $\a$ and $\s$ respectively. 
Any typical element $na = \exp (X+Z) a$ of $NA$ is denoted by $(X,Z,a)$
and the product law in $NA$ is given by
$$(X,Z,a)(X',Z',a')=(X+a^{1/2}X',
Z+aZ'+\frac{1}{2}a^{1/2}[X,X'],aa').$$  One can endow $NA$ with a
suitable left-invariant Riemannian metric which makes it a harmonic
manifold \cite{DR2}. Via the map $$h(X,Z,a) = (X,Z,
a+\frac{1}{4}\mid X\mid^2)$$ we can identify $S$ with the
generalized Siegel domain ${\mathcal D} = \{(X,Z,a) \in {\s} : a >
\frac{1}{4}\mid X\mid^2 \}$. Under this identification, $N$ gets
identified with the boundary $\partial {\mathcal D}$ of ${\mathcal
D}$.

Let ${\mathcal L}$ be the Laplace-Beltrami operator on $S$ with
respect to the Riemannian structure on $S$. Then, by Theorem 2.1 of \cite{D}, 
${\mathcal
L}= \sum\limits_{i=1}^{2m+k} E_i^2+E_o^2 -Q E_o$ where $E_1,
\ldots, E_{2m}$ in $\v$, $E_{2m+1},\ldots, E_{2m+k}$ in $\z$, $E_o$
in $\a$ form an orthonormal basis of $\s =\v \oplus \z \oplus \a$.
Bounded harmonic functions $u$ on ${\mathcal D}$ are those
functions that satisfy ${\mathcal L} u=0$ and have boundary values
almost everywhere on $\partial {\mathcal D}$ that is,
\begin{align}\label{e3.4}
\lim\limits_{a \rightarrow 0} u(X,Z,a) = \phi(X,Z)~~\mbox{a.e.},
\end{align}
where $\phi \in L^\infty (N)$(-see Theorem 3.7 of \cite{D2}). Moreover
$$u(X,Z,a) = (\phi \ast P_a)(X,Z),~~(X,Z,a) \in S,$$ where $P_a$ is
the Poisson kernel on the nilpotent group $N$ given by
$$P_a(X,Z) = \displaystyle{\frac{Ca^Q}{\lb \lb a+\frac{\mid
X\mid^2}{4}\rb^2 + \mid Z\mid^2\rb ^{Q}}} =
\displaystyle{\frac{Ca^Q}{{\lb \lb a+\frac{\mid X\mid^2}{4}\rb^2 +
\mid Z\mid^2\rb ^{m+k}}}}$$ and the convolution is on $N$. Here the constant $C$ is
chosen in 
such a way that $\|P_a\|_1 =1$. Note that
$P_a \in L^1(N)^\natural$ and
$$P_a(X,Z) = \displaystyle{\frac{Ca^Q}{{\lb a+\frac{\mid
X\mid^2}{4}+i \mid Z\mid\rb^{m+k} }\lb a+\frac{\mid X\mid^2}{4} -
i \mid Z\mid\rb^{m+k}}}.$$ \vskip0.15in

For each $a \in A$, we know that $\delta_a$ is an automorphism of
$N$, hence $\delta_a$ defines an automorphism of $L^1(N)^\natural$
by
$$ (\delta_a f)(X,Z) = a^{m+k} f(\delta_a(X,Z)).$$
Let $m(L^1(N)^\natural)$ be the set of non-zero
multiplicative linear functionals on $L^1(N)^\natural$. Then
$\delta_a$ induces a map $\delta_a^{\ast}$ on $m(L^1(N)^\natural)$ by
$$<f, \delta_a^{\ast} \psi> = <\delta_a f, \psi>,~~~\psi \in
m(L^1(N)^\natural),~~f \in L^1(N)^\natural.$$ Easy to see
that $\delta_a^{\ast}$ maps $m(L^1(N)^\natural)$
homeomorphically onto itself. If $f \in L^1(N)$ and $\|f\|_1 = 1$ then we can see
that $\{\delta_a f\}$ is an approximate identity in $L^1(N)$ as $a\rightarrow 0$.

We shall now make a small computation which we need in the next
section (see Lemma 3.4, \cite{CH}). For $w \in S_\z$, we shall
denote by $w^{\bot}$ the orthogonal compliment of $w$ in $\z$. Then
\begin{align}\label{e3.1}&\int\limits_{\z} e^{-i \nu \lan Z,
w\ran_{\n}}P_a(X,Z) \,dZ\\ \nonumber &= \int\limits_{\exp(w^\bot)}
\int\limits_{\bR}e^{-i \nu \lan tw+Z', w\ran_{\n}}P_a(X,tw+Z')\, dt
dZ'\\ \nonumber &= \int\limits_{\bR} e^{-i \nu t}
\int\limits_{\exp(w^\bot)} P_a(X,tw+Z') \,dZ' dt\\ \nonumber &=
\int\limits_{\bR} e^{-i \nu t} \int\limits_{\exp(w^\bot)} Ca^Q \lb \lb a
+\displaystyle{\frac{\mid X\mid^2}{4}}\rb^2 + t^2 + \mid Z'\mid^2
\rb^{-m-k} \,dZ' dt\\ \nonumber &= \int\limits_{\bR} e^{-i \nu t}
Ca^Q \int\limits_{\exp(w^\bot)}\lb  u^2 \rb^{-m-k} \lb 1+ \lb
\displaystyle{\frac{\mid Z'\mid}{u}}\rb^2 \rb^{-m-k}\, dZ' dt
\end{align}
where $u^2 = \lb \lb a +\displaystyle{\frac{\mid X\mid^2}{4}}\rb^2
+ t^2 \rb$. Now by a change of variable argument, we can show that
the above integral is
\begin{align}\label{e3.2}
= \int\limits_{\bR} e^{-i \nu t} C' a^Q \lb \lb a
+\displaystyle{\frac{\mid X\mid^2}{4}}\rb^2 + t^2
\rb^{-\frac{(m+k+1)}{2}} dt
\end{align}

In the next section, we shall prove our main result as a consequence
of a Tauberian theorem on H-type groups. \vskip0.25in
\section{Main Result}
\vskip0.15in \hspace{1.5cm}
In this section, we shall show that $R = L^1(N)^\natural$ is a regular $*$-algebra,
state and prove a Wiener - Tauberian theorem for $L^1(N)$. Further, we show that the
Gelfand transform of the Poisson kernel never vanishes. We shall conclude our main
result as a consequence of the Wiener - Tauberian theorem. Our proof of the Wiener -
Tauberian theorem is based on that of Hulanicki-Ricci \cite{HR}.
\vskip0.15in
\begin{proposition}
The commutative Banach algebra $R$ is regular.
\end{proposition}
\vskip0.15in
\noindent
{\bf Proof:} Given $f \in R$, define $f^*$ by $f^*(n) = \overline{f(n^{-1})}$. It is
easy to see that $f^* \in R$, $*$ defines an involution on $R$ and $R$ is symmetric.
Let $\tilde{R}$ be the commutative $*$-Banach algebra obtained from $R$ by adjoining
the unit element $1$. As in \cite{F2}, we can check that the set of multiplicative
linear functionals $m(\tilde{R})$ on $\tilde{R}$ is actually equal to $m(R) \cup \{
\infty \}$. Note that $m(\tilde{R})$ is compact and $\hat{R}$ separates points.
Further, if $f \in R \subset \tilde{R}$, then $\hat{f}(\infty) = 0$. Since $R$ is
$*$-closed, $\hat{R}$ is self-adjoint.

Let $C \subset m(R)$ be closed and $\xi \in m(R)\setminus C$. To show that $R$ is
regular, we need to show that there exists $f \in R$ such that $\hat{f}(C) = 0$ but
$\hat{f}(\xi) =1$. Since $C$ is closed in $m(R)$, $C \cup \{ \infty \}$ is compact
in $m(\tilde{R})$ and $\xi \notin C \cup \{ \infty \}$. As $m(\tilde{R})$ is compact
and Hausdorff, by Urysohn's lemma, we can get a continuous function $\phi$ on
$m(\tilde{R})$ such that $\phi(C) = 0, \phi(\infty) =0$, but $\phi(\xi) = 1$. By
Theorem 2 in \cite{N} (pp.217), closure of $\tilde{R} = C_{\infty}(m(\tilde{R})) =
C_0(m(R))$. So there exists $f \in R$ such that $\sup\limits_{\eta \in m(\tilde{R})}
\mid \phi(\eta) - \hat{f}(\eta) \mid < \displaystyle{\frac{1}{4}}$. Let $F$ be a
smooth real valued function defined on $\bR$ as follows:\\ \begin{align*}
F(\alpha) &= 0~~~\mbox{if}~~~\alpha \leq \frac{1}{3},\\ &=
1~~~\mbox{if}~~~\frac{2}{3} \leq \alpha \leq \frac{4}{3}.\end{align*} Then by
Dixmier (\cite{Di}), $F\circ \hat{f} \in \hat{R}$. But if $\eta \in C$, then $\mid
\hat{f}(\eta) \mid < \frac{1}{4} < \frac{1}{3}$ and so $(F \circ \hat{f})(\eta) =0$
and $\mid \hat{f}(\xi) - 1\mid < \frac{1}{4}$. This implies that $\frac{3}{4} <
\hat{f}(\xi) < \frac{5}{4}$ which in turn implies that $(F \circ \hat{f})(\xi) = 1$.
This proves the proposition.

For $f \in R$, let $C_f = \mbox{Support of }\hat{f}$ in $m(R)$. Let ${\mathscr B}
=$\linebreak $\lg f \in R : C_f \mbox{ is compact in }m(R) \rg$. Now we shall state
and prove a Wiener - Tauberian theorem for $R$.
\vskip0.15in
\begin{proposition}\label{p4.2}
The set ${\mathscr B}$ is dense in $R$.
\end{proposition}
\vskip0.15in
\noindent
{\bf Proof:} Let $f$ be a radial function in $L^1(N)$ having compact support. Choose
a function $F$ in $C^k(\bR)$ such that $F(1)=1$ and $F(x) = 0$ for $\mid x\mid \leq
\frac{1}{2}$. By \cite{Di}, we know that $\hat{f_1} = F \circ \hat{f} \in \hat{R}$
for some $f_1 \in R$. Note that $\int\limits_N f_1(X,Z) dX dZ = (F\circ
\hat{f})(0,0) = 1$ and $C_{f_1}$ is compact. Therefore $\widehat{(\delta_a f_1)}$
has compact support in $m(R)$ for every $a >0$. We know that $\{ \delta_a f_1\}$ is
an approximate identity in $L^1(N)$ as $a \rightarrow 0$. Therefore $f \ast \delta_a
f_1 \rightarrow f$ as $a \rightarrow 0$. But $\widehat{(f\ast \delta_a f_1)} =
\hat{f} \widehat{(\delta_a f_1)}$ has compact support in $m(R)$. This proves that
${\mathscr B}$ is dense in $C_c(N)$. But $C_c(N)$ is dense in $L^1(N)$, hence
${\mathscr B}$ is dense in $L^1(N)$.

As a consequence of the above we have the following:
\vskip0.15in
\begin{proposition}\label{p4.3}
Let $I$ be a proper closed right ideal in $L^1(N)$. Then there exists a $\psi \in
m(R)$ such that for all $f \in I\cap R$, $\hat{f}(\psi) = \psi(f) = 0$.
\end{proposition}
\vskip0.15in
\noindent
{\bf Proof:} Since $I$ is a proper closed right ideal in $L^1(N)$, an approximate
identity argument shows that $I\cap R$ is a proper closed ideal in $R$. 

To prove the proposition we need the following local Wiener- Tauberian theorem (-see
\cite{M}): Suppose $G \in R$ with $C_G$ compact. Let $f \in R$ be such that
$\hat{f}(\psi) \neq 0$ for all $\psi \in C_G$. Then there exists $g \in R$ such that
$g \ast f \ast G =G$. In order to prove this claim, note that $\hat{f} $ is
continuous on the compact set $C_G$. Hence there exists $\beta > 0$ such that
$\hat{f}(\phi) > \beta$ for all $\phi \in C_G$. Now choose $F$ in $C^k(\bR)$ such
that 
\begin{align*}
F(\alpha)&= \displaystyle{\frac{1}{\alpha}}~~~\mbox{for}~~~\alpha \geq \beta,\\&=
0~~~~~\mbox{for}~~~~ \alpha \leq 0.
\end{align*}
Then $F \circ \hat{f}|_{C_G} = \displaystyle{\frac{1}{\hat{f}}}$. But by \cite{Di},
we have $g \in R$ such that $\hat{g} = F \circ \hat{f} \in \hat{R}$. Therefore
\begin{align*}\widehat{g\ast f \ast G}(\psi) &= (\hat{g}\hat{f}\hat{G})(\psi)\\
&=\begin{cases} 0~~~~\mbox{if}~~~~ \psi \notin C_G,\\ \hat{G}~~~~\mbox{if}~~~~\psi
\in C_G\end{cases}\end{align*}
that is, $\hat{g}\hat{f}\hat{G}=\hat{G}$. By uniqueness of Gelfand transform $g \ast f
\ast G =G$.

We shall now prove the proposition. Assume on the contrary that for every $\psi \in
m(R)$, there exists $f$ in $I \cap R$ such that $\hat{f} (\psi) \neq 0$. Let $G \in
R$ be such that $C_G$ is compact. By assumption, given any $\psi \in C_G$ we can
find $f_{\psi} \in I \cap R$ such that $\hat{f_{\psi}}(\psi) \neq 0$. In fact
$\hat{f_{\psi}}(\psi) > 0$. By continuity, $\hat{f}$ does not vanish in a
neighborhood $U_{\psi}$ of $\psi$. The collection of open sets$\{U_{\psi}\}_{\psi
\in C_G}$ forms an open cover for $C_G$. Hence we can find $\psi_1, \psi_2, \cdots,
\psi_n$ in $C_G$ such that $\{U_{\psi_1}, U_{\psi_2}, \cdots, U_{\psi_n} \}$ forms a
finite subcover of $C_G$. Consider the function $f = f_{\psi_1} + f_{\psi_2} +
\cdots + f_{\psi_n}$. Then $f \in I \cap R$ as $I$ is an ideal and $\hat{f}(\psi)
\neq 0$ for all $\psi \in C_G$. By our claim above, then there exists $g \in R$ such
that $g \ast f \ast G =G$. But $I \cap R$ is an ideal in $R$ and $f \in I \cap R$.
Hence $g \ast f \ast G \in I \cap R$. This implies that $G\in I \cap R$. This in
turn implies that ${\mathscr B} \subset I \cap R$. Therefore, $R =
\overline{\mathscr B} \subset \overline{I \cap R} = I \cap R$ as $I \cap R$ is
closed. This contradicts the fact that $I \cap R$ is a proper closed ideal in $R$.
This completes the proof of the proposition.

Put $s = \frac{m+k-1}{2}$. Note that $s \geq 0$. For $a \in A$, recall that the Poisson kernel $P_a$ is given by 
\begin{align*}
P_a(X,Z) &= \displaystyle{\frac{Ca^Q}{\lb \lb a+\frac{\mid X\mid^2}{4}\rb^2 + \mid
Z \mid^2 \rb^{2s+1}}},~~(X,Z) \in N, \\ &=
\displaystyle{\frac{Ca^Q}{(2s!)^2}\frac{(2s!)^2}{\lb a+\frac{\mid X\mid^2}{4} -
i\mid Z\mid \rb^{2s+1}\lb a+\frac{\mid X\mid^2}{4} + i\mid Z\mid \rb^{2s+1}}}
\end{align*}
Consider for any $r>0$, $$F_a (X,Z) = \displaystyle{\frac{1}{r!} \frac{r!}{\lb
a+\frac{\mid X\mid^2}{4} + i\mid Z\mid \rb^{r+1}}}.$$ Using the Laplace transform
techniques, one can easily show that
$$F_a(X,Z) = \displaystyle{\frac{1
}{r!}\int\limits_0^{\infty} e^{i \alpha \mid Z\mid} e^{-\alpha (a+\frac{\mid
X\mid^2}{4}+ i \mid Z\mid)}\alpha^r \,d\alpha}.$$
Therefore,
\begin{align}\label{e4.6}
P_a(X,Z)&=
\displaystyle{\frac{Ca^Q}{(2s!)^2}\int\limits_0^{\infty}\int\limits_0^{\infty}
e^{-\alpha (a + \frac{\mid X\mid^2}{4} - i\mid Z\mid)} e^{-\beta (a + \frac{\mid
X\mid^2}{4} + i\mid Z\mid)} (\alpha \beta)^{2s} \,d\alpha d\beta }\nonumber\\ &=
\displaystyle{\frac{1}{(2s!)^2}\int\limits_0^{\infty}\int\limits_0^{\infty}
e^{i(\alpha-\beta) \mid Z\mid} e^{-(\alpha +\beta) (a + \frac{\mid X\mid^2}{4})}
(\alpha \beta)^{2s} \,d\alpha d\beta }
\end{align}
As the last step in our proof of the main result, we have the following:
\vskip0.15in
\begin{proposition}\label{p4.4}
For every $a >0$, the Gelfand transform $\hat{P_a}$ of $P_a$ is never zero on $m(R)$.
\end{proposition}
\vskip0.15in
\noindent
{\bf Proof:} We need to check that $\hat{P_a}(\nu, l)$ and $\hat{P_a}(0,\mu)$ do not
vanish for $\nu > 0$, $\mu \geq 0$ and $l \in \bN$. Consider the integral
\begin{align}\label{e4.7}
&\int\limits_{\z}P_a(X,Z) e^{-i \nu \lan Z,w\ran_{\n}}\,dZ \nonumber\\&=
\int\limits_{\bR} \int\limits_{\exp (w^\bot)} P_a(X,tw+Z^\prime) e^{-i \lan
tw+Z^\prime, \nu w\ran_\n} \,dZ^\prime dt\nonumber\\&= \int\limits_{\bR} e^{-i \nu t}
\int\limits_{\exp (w^\bot)} P_a (X,tw+Z^\prime) e^{-i \lan Z^\prime, \nu w\ran_\n}
\,dZ^\prime dt\nonumber\\ &= \int\limits_{\bR}e^{-i \nu t} \int\limits_{\exp
(w^\bot)} P_a (X,tw+Z^\prime) \,dZ^\prime dt \nonumber\\ &= Ca^Q
\int\limits_{\bR}e^{-i \nu t} \lb \lb a+\frac{\mid X\mid^2}{4} \rb^2 + t^2 \rb
^{-(s+1) }\,dt 
\end{align}
by (\ref{e3.1}). Consider the expression
\begin{align}\label{e4.8}
\lb \lb a +\frac{\mid X\mid^2}{4}\rb^2 + t^2 \rb^{-(s+1)} &= \lb  a + \frac{\mid
X\mid^2}{4} +it\rb^{-(s+1)}\lb  a + \frac{\mid X\mid^2}{4}
-it\rb^{-(s+1)}\nonumber\\ &= \frac{1}{(s!)^2}
\int\limits_0^{\infty}\int\limits_0^{\infty} e^{i(\alpha - \beta)t}
e^{-(\alpha+\beta)\lb a+ \frac{\mid X\mid^2}{4}\rb} \alpha^s \beta^s \,d\alpha d\beta
\end{align}
by equation (\ref{e4.6}). Now using Fourier transform techniques together with
(\ref{e4.7}), (\ref{e4.8}), we get,
\begin{align}\label{e4.9}
\int\limits_{\z}P_a(X,Z) e^{-i \nu \lan Z,w\ran_{\n}}dZ=\frac{Ca^Q}{(s!)^2}
\int\limits_0^{\infty} e^{-(2\beta+\nu)\lb a+ \frac{\mid X\mid^2}{4}\rb}\lb \beta
+\nu\rb^s \beta^s \,d\beta.
\end{align} 
Now if we take $\nu = 0$ and evaluate the Fourier transform in the variable $X$, we 
obtain $\hat{P_a}(0,\mu)$.

Therefore 
\begin{align}
\hat{P_a}(0,\mu) &= \hat{P_a}(0,\mid Y\mid)~~~~~~\mbox{~~where~~}Y \in \v\nonumber\\&= \frac{Ca^Q}{(s!)^2} \int\limits_{\v}\int\limits_0^\infty
e^{-2\beta(a+\frac{\mid X\mid^2}{4})} \beta^{2s} e^{-i\lan  Y, X\ran} \,d\beta
dX\nonumber\\&= 
\frac{Ca^Q}{(s!)^2}\int\limits_0^\infty e^{-2\beta a} \beta^{2s} \int\limits_{\v}
e^{-\beta \frac{\mid X\mid^2}{2}} e^{i \lan Y,X\ran } \,dX d\beta \nonumber\\ &= 
\frac{Ca^Q}{(s!)^2}\int\limits_0^\infty e^{-\frac{ \mu^2}{2\beta}} e^{-2\beta a}
\beta^{2s} \,d\beta 
\end{align}
Let $I = \int\limits_0^\infty e^{-2\beta a} \beta^{2s} e^{-\frac{\mu^2}{2\beta}}
\,d\beta$. Choose $0<\epsilon_1 <\epsilon_2 <\infty$. Then $I = $ \linebreak
$\int\limits_0^{\epsilon_1} + \int\limits_{\epsilon_1}^{\epsilon_2}
+\int\limits_{\epsilon_2}^\infty e^{-2\beta a} \beta^{2s} e^{-\frac{\mu^2}{2\beta}}
\,d\beta$. Note that for $\beta$ in $[\epsilon_1,\epsilon_2]$ the integrand
\linebreak
$e^{-2\beta a} \beta^{2s} e^{-\frac{\mu^2}{2\beta}} > 0$ and the integrand is
non-negative for all other values in the interval $[0,\infty)$. Hence
$\int\limits_{\epsilon_1}^{\epsilon_2}e^{-2\beta a} \beta^{2s}
e^{-\frac{\mu^2}{2\beta}}\,d\beta > 0$ and the integral on the intervals $[0, \epsilon_1)$,
$[\epsilon_2,\infty)$ are non-negative. Therefore $I > 0$ which in turn implies that
$\hat{P_a}(0,\mu)>0$. For $\nu \neq 0$, we have
\begin{align}\label{e4.11}
\hat{P_a}(\nu,l) &= \frac{Ca^Q}{(s!)^2} \int\limits_{\v} \int\limits_0^\infty
e^{-(2\beta + \nu)(a + \frac{\mid X\mid^2}{4})} (\beta +\nu)^s \beta^s  
e^{-\frac{\nu}{4} \mid X\mid^2} \displaystyle{\frac{L^{m-1}_l(\frac{\nu}{2}  \mid
X\mid^2)}{\lb \begin{smallmatrix}l+m-1\\l\end{smallmatrix}\rb} } \,d\beta
dX\nonumber\\ &= \frac{Ca^Q}{(s!)^2} \int\limits_0^\infty
e^{-(2\beta + \nu)a} (\beta +\nu)^s \beta^s \int\limits_{\v} e^{-\frac{\nu}{4} \mid
X\mid^2} \displaystyle{\frac{L^{m-1}_l(\frac{\nu}{2}  \mid X\mid^2)}{\lb
\begin{smallmatrix}l+m-1\\l\end{smallmatrix}\rb} }\,dX d\beta \nonumber\\ &=
\frac{Ca^Q}{(s!)^2} \frac{\mid S_{\v} \mid}{\lb
\begin{smallmatrix}l+m-1\\l\end{smallmatrix}\rb} \int\limits_0^\infty e^{-(2\beta +
\nu)a} (\beta +\nu)^s \beta^s \lb \int\limits_0^\infty e^{-\frac{\nu
k^2}{4}}L^{m-1}_l\lb \frac{\nu k^2}{2} \rb k^{2m-1} \,dk \rb d\beta \nonumber\\ &=
\frac{Ca^Q}{(s!)^2} \frac{\mid S_{\v} \mid}{\lb
\begin{smallmatrix}l+m-1\\l\end{smallmatrix}\rb}  \frac{2^{m-1}}{\nu^m}
\int\limits_0^\infty e^{-(2\beta + \nu)a} (\beta +\nu)^s \beta^s
\int\limits_0^\infty e^{-\frac{y}{2}} L^{m-1}_l(y) y^{m-1} \,dy d\beta
\end{align}
by a change of variable. But we know from \cite{T}, that 
\begin{align*}
L^{m-1}_l (x) = \frac{1}{l!} e^x x^{-(m-1)} \frac{d^l}{dx^l}\lb e^{-x} x^{l+m-1} \rb 
\end{align*}
Integrating by parts (\ref{e4.11}) implies
\begin{align}
\hat{P_a}(\nu,l) &= \frac{Ca^Q}{(s!)^2} \frac{\mid S_{\v} \mid}{\lb
\begin{smallmatrix}l+m-1\\l\end{smallmatrix}\rb}  \frac{2^{m-1}}{\nu^m}
\int\limits_0^\infty e^{-(2\beta + \nu)a} (\beta +\nu)^s \beta^s
\nonumber\\&\int\limits_0^\infty e^{-\frac{y}{2}} e^y y^{-(m-1)}
\frac{d^l}{dy^l}\lb e^{-y} y^{l+m-1} \rb y^{m-1} \,dy d\beta \nonumber\\ &=
\frac{Ca^Q}{(s!)^2} \frac{\mid S_{\v} \mid}{\lb
\begin{smallmatrix}l+m-1\\l\end{smallmatrix}\rb}  \frac{2^{m-1}}{\nu^m}
\int\limits_0^\infty e^{-(2\beta + \nu)a} (\beta +\nu)^s \beta^s
\int\limits_0^\infty e^{\frac{y}{2}} \frac{d^l}{dy^l}\lb e^{-y} y^{l+m-1} \rb \,dy
d\beta \nonumber\\ &= \frac{Ca^Q}{(s!)^2} \frac{\mid S_{\v} \mid}{\lb
\begin{smallmatrix}l+m-1\\l\end{smallmatrix}\rb}  \frac{2^{m-1}}{\nu^m} \lb
\frac{-1}{2} \rb^l \int\limits_0^\infty e^{-(2\beta + \nu)a} (\beta +\nu)^s \beta^s
\int\limits_0^\infty e^{\frac{y}{2}} e^{-y} y^{l+m-1} \,dy d\beta \nonumber\\ &=
\frac{Ca^Q}{(s!)^2} \frac{\mid S_{\v} \mid}{\lb
\begin{smallmatrix}l+m-1\\l\end{smallmatrix}\rb}  \frac{2^{m-1}}{\nu^m} \lb
\frac{-1}{2} \rb^l \int\limits_0^\infty e^{-(2\beta + \nu)a} (\beta +\nu)^s \beta^s
\int\limits_0^\infty e^{\frac{-y}{2}} y^{l+m-1} \,dy d\beta \nonumber\\ &=
\frac{Ca^Q}{(s!)^2} \frac{\mid S_{\v} \mid}{\lb
\begin{smallmatrix}l+m-1\\l\end{smallmatrix}\rb} \frac{2^{2m-1}(-1)^l}{\nu^m (l!)}
(l+m-1)! \int\limits_0^\infty e^{-(2\beta +\nu) a} (\beta+\nu)^s \beta^s \,d\beta.
\end{align}
By repeating a similar argument as in the case of $\hat{P_a}(0,\mu)$ we can show
that $\hat{P_a}(\nu,l) > 0$. This completes the proof of the assertion.

We shall prove the tangential convergence of the bounded harmonic functions on the
generalised Siegel domain ${\mathcal D}$.
\vskip0.15in
\begin{theorem}
Suppose $u$ is a bounded harmonic function on ${\mathcal D}$ and
$$\lim\limits_{(X,Z) \rightarrow \infty} u(X,Z,a_o) = \alpha$$ for some $a_o > 0$.
Then for all $a > 0$, the limit $\lim\limits_{(X,Z) \rightarrow \infty} u(X,Z,a)$
exists and is equal to $\alpha$. 
\end{theorem}
\vskip0.15in
\noindent
{\bf Proof:} Let $\phi$ be a function in $L^\infty(N)$ satisfying (\ref{e3.4}).
Consider the right ideal in $L^1(N)$ given by $$I = \lg g \in L^1(N) : \phi \ast g =
\alpha \int\limits_N g(X,Z) dX dZ \rg.$$ By our assumption, $P_{a_o} \in I\cap R$.
But by Proposition \ref{p4.4}, $\hat{P_{a_o}}$ does not vanish anywhere on the
Gelfand spectrum $m(R)$ of $R$. By Proposition \ref{p4.3}, this would imply that $I
= L^1(N)$. This completes the proof of the theorem.

\end{document}